\documentclass[a4paper%,draft
]
{book}
\usepackage[T2A,TS1]{fontenc}
\usepackage[cp1251]{inputenc}
\usepackage[english,russian]{babel}
\usepackage{amsfonts,amsmath,amssymb,amsthm,graphicx,afterpage}
\makeatletter

\def\hm#1{#1\nobreak\discretionary{}{\hbox{$#1$}}{}}

\let\Bbb=\mathbb
\long\def\comment#1\endcomment{}
\def\R{{\mathbb R}} \def\Z{{\mathbb Z}}  \def\Q{\Bbb Q}
\newcommand{\diam}{\mathop{\mathrm{diam}}}

\newskip\thmindent
\newtheoremstyle{myplain}% name
  {\thmindent}%      Space above
  {\thmindent}%      Space below
  {\itshape}%         Body font
  {\parindent}%         Indent amount (empty = no indent, \parindent = para indent)
  {\bfseries}% Thm head font
  {.}%        Punctuation after thm head
  { }%     Space after thm head: " " = normal interword space;
  {}%         Thm head spec (can be left empty, meaning `normal')

\newtheoremstyle{mydefinition}% name
  {\thmindent}%      Space above
  {\thmindent}%      Space below
  {\normalfont}%         Body font
  {\parindent}%         Indent amount (empty = no indent, \parindent = para indent)
  {\bfseries}% Thm head font
  {.}%        Punctuation after thm head
  { }%     Space after thm head: " " = normal interword space;
  {}%         Thm head spec (can be left empty, meaning `normal')

\theoremstyle{myplain}
\newtheorem{theorem}{Теорема}
\newtheorem*{theorem*}{Теорема}
\newtheorem*{theoremB}{Теорема Бэра о категории}

\newtheorem{statement}[theorem]{Утверждение}
\newtheorem{corollary}[theorem]{Следствие}
\newtheorem{conjecture}[theorem]{Гипотеза}

\theoremstyle{mydefinition}

\newtheorem*{definition*}{Определение}
\newcommand{\qdefinitions}{\indent\textbf{Определение.} \ignorespaces}
\newtheorem{remark}[theorem]{Замечание}
\newtheorem{problem}{}

\renewenvironment{proof}[1][Доказательство]{\par%\addvspace{\medskipamount}
  \pushQED{\qed}%
  \normalfont %\topsep6\p@\@plus6\p@\relax
\indent\textit{#1.} \ignorespaces
}{%
  \popQED%\endtrivlist\@endpefalse
  \par%\addvspace{\medskipamount}
}

\renewcommand\section{\@startsection {section}{1}{\parindent}%
                                   {2.5ex \@plus 1ex \@minus .2ex}%
                                   {1.5ex \@plus.2ex}%
                                   {\centering\normalfont\large\bfseries}}
\renewcommand\subsection{\@startsection{subsection}{2}{\parindent}%
                                     {1.5ex\@plus .5ex \@minus .2ex}%
                                     {1ex \@plus .2ex}%
                                     {\centering\normalfont\bfseries}}
\newcommand{\tsubsection}[1]{\subsection*{\hskip-\leftmargin #1}}

\newcounter{myitemize}
\newenvironment{myitemize}{\nopagebreak\setcounter{myitemize}{0}\def\item{\par\refstepcounter{myitemize}%
(\alph{myitemize}) \ignorespaces}}{}

\let\emptyset\varnothing
\let\dots\ldots 
\let\ge\geqslant
\let\le\leqslant

\newcommand{\?}{\nobreak\hskip.167em\nobreak\hskip\z@skip}

\newcommand{\hsk}[1]{\medmuskip4mu minus #1mu
\thickmuskip4mu minus #1mu}

\renewenvironment{thebibliography}[1]
     {\section*{\bibname}%
      \addcontentsline{toc}{section}{\bibname}
      \list{\@biblabel{\@arabic\c@enumiv}}%
           {\settowidth\labelwidth{\@biblabel{#1}}%
            \leftmargin\labelwidth
            \advance\leftmargin\labelsep
            \@openbib@code
            \usecounter{enumiv}%
            \let\p@enumiv\@empty
            \renewcommand\theenumiv{\@arabic\c@enumiv}}%
      \sloppy
      \clubpenalty4000
      \@clubpenalty \clubpenalty
      \widowpenalty4000%
      \sfcode`\.\@m}
     {\def\@noitemerr
       {\@latex@warning{Empty `thebibliography' environment}}%
      \endlist}

\newcommand{\textsum}{\sum\limits}
\let\epsilon\varepsilon

\long\def\@makecaption#1#2{%
  \vskip\abovecaptionskip
  \sbox\@tempboxa{\small #1. #2}%
  \ifdim \wd\@tempboxa >\hsize
    \unhbox\@tempboxa\par
  \else
    \global \@minipagefalse
    \hbox to\hsize{\hfil\box\@tempboxa\hfil}%
  \fi
  \vskip\belowcaptionskip}

\newlength{\wrfwidth}
\newlength{\tempa}
\newlength{\tempb}
\newlength{\tempc}

\newcommand{\from}{/\!/ }

\renewcommand\tableofcontents{%
    \if@twocolumn
      \@restonecoltrue\onecolumn
    \else
      \@restonecolfalse
    \fi
    \section*{\contentsname}
    \@starttoc{toc}%
    \if@restonecol\twocolumn\fi
    }
\renewcommand\normalsize{%
   \@setfontsize\normalsize\@xpt\@xiipt
   \abovedisplayskip 6\p@ \@plus3\p@ \@minus1.5\p@
   \abovedisplayshortskip \z@ \@plus3\p@
   \belowdisplayskip\abovedisplayskip
   \belowdisplayshortskip\belowdisplayskip
   \let\@listi\@listI
}
\normalsize

\def\headbreak{\\}

\renewcommand*\l@section{\medskip\@dottedtocline{1}{1.5em}{2.3em}}

\let\mytexttt\relax

\declare@shorthand{russian}{"=}{\nobreak\hskip\z@skip\hbox{-}\nobreak\hskip\z@skip}

\def\ps@plain{\let\@mkboth\@gobbletwo
     \let\@oddhead\@empty\def\@oddfoot{\reset@font\small\hfil\thepage
     \hfil}\let\@evenhead\@empty\let\@evenfoot\@oddfoot}
\advance\footskip-2pt
\makeatother

\usepackage[perpage]{footmisc}
\usepackage{wrapfig}

\begin{document}

%Ambient homogeneity
%This note is purely expository. A subset $N$ of the plane is affine ambient
%homogeneous if for each $x,y\in N$ there exists an affine transformation taking
%$x$ to $y$ and $N$ to itself. The result of D.\?Repovs, E.\?V.\?Scepin and the
%author on such subsets is presented, together with discussion, corollaries and
%generalizations. At the end some non"=elementary corollaries are presented
%(including a simple proof of the smooth version of the Hilbert"--~Smith conjecture
%on topological groups). Most part of the text is accessible to undergraduates
%familiar with the notion of continuity. The text could be an interesting easy
%reading for mature mathematicians.
%57R50; Secondary: 53A04, 54H11, 58A05

%\title{Объемлемая однородность}
%\author{А.\?Скопенков\footnote{skopenko@mccme.ru;
%http://dfgm.math.msu.su/people/skopenkov/PAPERSCI.pdf.
%Moscow State University, Independent University of Moscow and
%Moscow Institute of Open Education.  }}
%
%\maketitle

%\subsection*{Аннотация (для обложки)}
 
%%%{

%%%\thispagestyle{empty}

%%%\vbox to\textheight{
%\vskip-\headheight
%\vskip-\headsep

{\centering

Летняя школа <<Современная математика>>\\
Дубна, июль 2009

\vspace*{2cm}

{\large А.\,Б.\,Скопенков}

\vspace*{3cm}

{\huge Объемлемая однородность\par}

\vspace*{2cm}

}

Брошюра написана по материалам миникурса в летней школе <<Современная математика>> в Дубне
в 2009~г. и доклада на семинаре по геометрии им. И.\?Ф.\?Шарыгина в 2010~г.

Понятие объемлемой однородности возникает из простых <<физических>> вопросов.
%Однако их важность для математики раскрывется далее.
Введение доступно школьнику (кроме его последнего пункта, где требуется понятие непрерывного отображения между подмножествами плоскости).
Далее практически <<школьными>> методами  мы получим характеризацию объемлемо однородных подмножеств плоскости.
В этой части уже необходимо знакомство с открытыми и замкнутыми множествами на прямой и плоскости.
Затем выясняется, что понятие объемлемой однородности связано со многими
важными теориями и результатами "--- теорией динамических систем, многообразий и групп Ли, пятой проблемой Гильберта и  проблемой Гильберта--Смита.
Приложение доступно студенту, знакомому с этими понятиями.

Брошюра адресована широкому кругу людей, интересующихся математикой.
Она может быть интересным <<легким чтением>> для профессиональных математиков.

\vfill

{\centering

Москва\\
Издательство \hbox{МЦНМО}\\
2012

}

%%%\thispagestyle{empty}

%%%}

%%%}

\vspace*{3cm}

%%%\clearpage\thispagestyle{empty}

\comment

\def\udk{???.??}
\def\bbk{??.??}
\def\authorscode{С44}
\def\ISBN{978-5-94057-905-2}

\thispagestyle{empty}

\vbox to\textheight{
%\vspace{-\headsep} \vspace{-.7\headheight}
\renewcommand{\bfdefault}{b}

\noindent\makebox[\textwidth][s]{\begin{tabular}[t]{@{}l@{ }l@{}}
УДК&\udk\\
ББК&\bbk\\
&\authorscode
\end{tabular}
\hfil \small }

\vspace{2\baselineskip}

{\centering

%\vspace{3\baselineskip}

\newdimen\kartwidth
\kartwidth\ifdim\textwidth>110mm 103mm\else\textwidth\fi
\fontsize{10}{10.5pt}\selectfont
\newdimen\acwidth\acwidth9mm\advance\acwidth5pt
\advance\kartwidth-\acwidth

\begin{tabular}{@{}p{\acwidth}@{}p{\kartwidth}@{}}
\raisebox{-11pt}[0pt][0pt]{\authorscode} &\parindent15pt

\medskip

\noindent\textbf{Скопенков~А.\,Б.}

Объемлемая однородность. "--- М.: \hbox{МЦНМО},
2012. "--- 32~с.

\vspace{4pt}

ISBN \ISBN

\vspace{4pt}

{\fontsize{8}{8.5pt}\selectfont
Брошюра написана по материалам миникурса в летней школе <<Современная математика>> в Дубне
в 2009~г. и доклада на семинаре по геометрии им. И.\?Ф.\?Шарыгина в 2010~г.

Понятие объемлемой однородности возникает из простых <<физических>> вопросов.
%Однако их важность для математики раскрывется далее.
Введение доступно школьнику (кроме его последнего пункта, где требуется понятие непрерывного отображения между подмножествами плоскости).
Далее практически <<школьными>> методами  мы получим характеризацию объемлемо однородных подмножеств плоскости.
В этой части уже необходимо знакомство с открытыми и замкнутыми множествами на прямой и плоскости.
Затем выясняется, что понятие объемлемой однородности связано со многими
важными теориями и результатами "--- теорией динамических систем, многообразий и групп Ли, пятой проблемой Гильберта и  проблемой Гильберта--Смита.
Приложение доступно студенту, знакомому с этими понятиями.

Брошюра адресована широкому кругу людей, интересующихся математикой.
Она может быть интересным <<легким чтением>> для профессиональных математиков.

% В частности, в приложении приводится простое доказательство гладкой версии гипотезы
% Гильберта--Смита о топологических группах, возникшей при решении 5-й проблемы Гильберта.
\par}

\rightline{ББК \bbk}
\end{tabular}

}

\vspace{\baselineskip}

\vspace{6pt}

{\small

}

\vfil

\vskip.5\baselineskip
\noindent\makebox[\textwidth][s]{%
{\begin{tabular}[b]{@{}l@{}}
%\includegraphics{174.eps}
%\small
\textbf{ISBN \ISBN}
\end{tabular}\hfil
\begin{tabular}[b]{@{}l@{ }l@{}}
\copyright& Скопенков~А.\,Б., 2012.\\
\copyright& \hbox{МЦНМО}, 2012.
\end{tabular}}}}

\endcomment

\clearpage

\hfill{\itshape Посвящается памяти В.\?И.\?Арнольда}

\vskip\baselineskip

\subsection*{Советы читателю}
Начать читать брошюру разумно с введения.
Его три пункта практически независимы друг от друга, и их можно читать в произвольном порядке.
Впрочем, они расположены в порядке возрастания сложности.
В дальнейшем из введения используется только пункт~\ref{subsec1.2}.

Оставшиеся параграфы практически независимы друг от друга, и их можно читать в произвольном порядке.
Впрочем, они расположены в порядке возрастания сложности.

Основное содержание брошюры "--- утверждение~\ref{p:diff-f} из пункта \ref{subsec1.2}, его доказательство в \ref{s:proof-aff} и его обобщения в параграфах \ref{sec4} и \ref{pril}.

В брошюре много задач, обозначаемых жирными цифрами.
Большинство задач несложны.
% ??? нужна ли первая запятая?
% Р
При этом, если условие задачи является формулировкой
утверждения, то это утверждение и надо доказать.
Формулировки задач нужно прочитать "--- это поможет вам
понять текст, даже если вы не сможете решить задачи.
Если некоторые встречающиеся, но не определенные понятия вам незнакомы, то можно или игнорировать соответствующую задачу,
%(приняв её утверждение на веру),
%(например, для понимания курса не нужно понятие действия группы на множестве)
или узнать определение (у преподавателя, в wikipedia, в книгах...).
% ??? заменить "мне" на "автору" (ср. ниже в разделе "Благодарности"?
% нет
Двумя звездочками отмечены задачи, решение которых мне неизвестно.

Обновляемая версия поддерживается на \mytexttt{http://arxiv.org/abs/1003.5278}.

\subsection*{Благодарности}

Автор благодарен В.\?Клепцыну, Г.\?Мерзону и А.\?Сосинскому за полезные замечания и обсуждения.
Автор был поддержан грантом фонда Саймонса.
%а также В.\?Клепцыну за подготовку рисунков.
% "--- с единственной оговоркой о том, какие предварительные сведения необходимы для их прочтения.
%Следующий раздел посвящён канторову множеству, после чего в разделе~\ref{s:proof-aff} мы докажем %теорему~\ref{thm:t1}. При этом нам потребуется утверждение о дифференцируемости липшицевой %функции, которое следует из теоремы Бэра о категории и которое мы разберём вместе с самой этой %теоремой.

\section{Введение}\label{sec1}

\subsection{Изометрическая объемлемая однородность}

\settowidth{\tempa}{\includegraphics{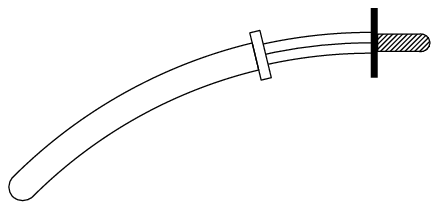}}%
\begin{wrapfigure}{o}{\tempa}
\vspace{-10pt}
\includegraphics{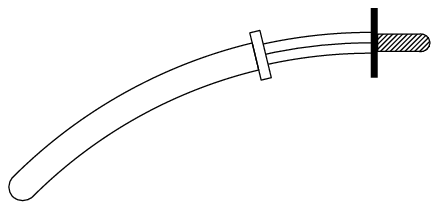}
\end{wrapfigure}
Какой формы могут быть ножны, чтобы из них можно было
вытащить саблю? Переформулируя этот вопрос на
математическом языке, мы приходим к следующему определению.
%\begin{figure}[h]\centering
%\includegraphics{sab.1}
%\end{figure}

%\begin{definition*}\label{def:isom}\end{definition*}

%\smallskip
%{\bf Предложение А.\?Б.\?Сосинского. Было бы замечательно
%поместить на обложке или здесь цветную фотографию сабли,
%частично вынутой из ножен.}
%(скачав из Интернета или из `Кванта'?).Если нельзя на обложке, то здесь. }

%\smallskip
%{\bf Определение.}
\begin{definition*}
Подмножество $N$ пространства $\R^m$ (в частности,
плоскости $\R^2$ или трехмерного пространства~$\R^3$)
называется \emph{изометрически объемлемо однородным}, если
для любых двух точек $x,y\in N$ существует движение (т.\,е.~изометрия)
пространства, переводящее $x$ в $y$, а $N$ в
себя.
%(Это движение не предполагается непрерывно зависящим от $x$ и $y$.)
\end{definition*}

\begin{figure}[t]
\vskip\topskip\vskip-6pt
\settowidth{\tempa}{\includegraphics[scale=.8]{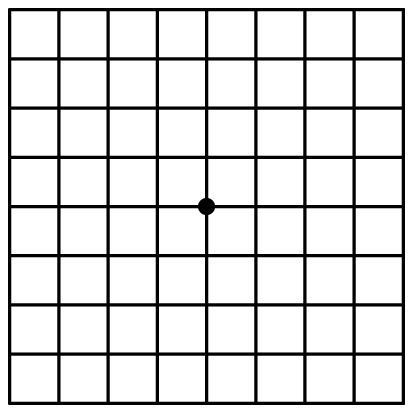}}
\settowidth{\tempb}{\small Рис. 2. Образ решетки при движении}
\hbox to\textwidth{\hss
\begin{minipage}[b]{\tempa}\centering
\includegraphics[scale=.8]{transform.1.eps}
\caption{Решетка
%и ее образ при движении, аффинном преобразовании и диффеоморфизме
}\label{fig:lattice}
\end{minipage}\hss
\begin{minipage}[b]{\tempb}\centering
\includegraphics[scale=.8]{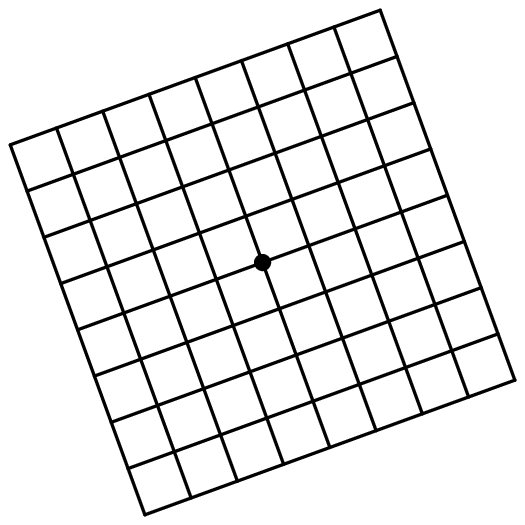}
\caption{\smash[b]{Образ решетки при движении}}\label{fig:latiso}
\end{minipage}%
\hss}
\end{figure}

%\begin{figure}[h]
%\hfill \includegraphics[scale=.5]{transform.2.eps}\quad
%\hfill \mbox{}
%\caption{Образ решетки при движении}\label{fig:latiso}
%\end{figure}

%\smallskip
Напомним, что движением (т.\,е.~изометрией) называется преобразование, сохраняющее расстояния,
см. рис. \ref{fig:lattice} и \ref{fig:latiso}.

Отметим, что в этом
%(традиционном)
определении не требуется непрерывной зависимости
%такого преобразования
движения от $x$ и $y$.
Хотя её и было бы естественно потребовать, исходя из исходной <<физической>> задачи.  %(определение~\ref{def:isom}, таким образом, лишь мотивировано исходным физическим вопросом).

%Приведём несколько примеров изометрически объемлемо однородных подмножеств.

%\smallskip
%{\it Рисунок: сфера, винтовая линия [Ar84, рис. 121] и тор в $\R^4$ [Ar84, рис. 171]}
%\smallskip

%Существенная часть брошюры, в особенности, разделы~
%\ref{s:Cantor} и~ \ref{s:Baire}, изложена в виде последовательности задач
% (третий и пятый пункты состоят из задач).
%{\it Общее замечание о задачах в этом тексте.}

\begin{figure}[t]\centering
\vskip\topskip\vskip-6pt
\settowidth{\tempa}{\includegraphics{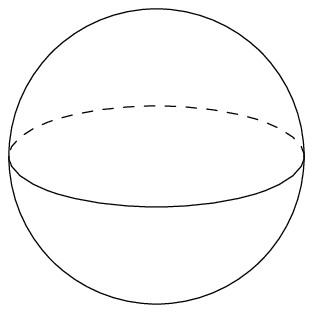}}
\settowidth{\tempb}{\includegraphics{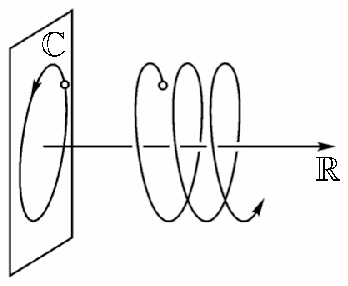}}
\settowidth{\tempc}{\includegraphics{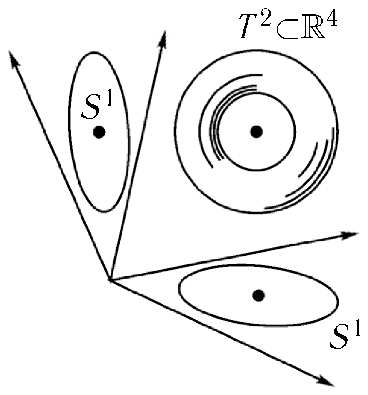}}
\hbox to\textwidth{\begin{minipage}[b]{\tempa}\centering
\includegraphics{sf.1.eps}
\caption{\smash[b]{Сфера}}\label{fig:sphere}
\end{minipage}\hss
\begin{minipage}[b]{\tempb}\centering
\includegraphics{121-new-tak-good.eps}\\
\caption{Винтовая линия}\label{fig:vint}
\end{minipage}\hss
\begin{minipage}[b]{\tempc}\centering
\includegraphics{171-new-corrected-good.eps}
\caption{\smash[b]{Тор в $\R^4$}}\label{fig:torR4}
\end{minipage}%
}
%\vskip-1pt
\end{figure}

\begin{figure}[t]\centering
\settowidth{\tempa}{\includegraphics{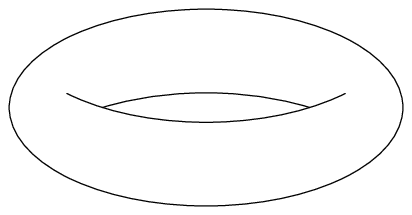}}
\hbox to\textwidth{\hss\begin{minipage}[b]{.49\textwidth}\centering
\includegraphics{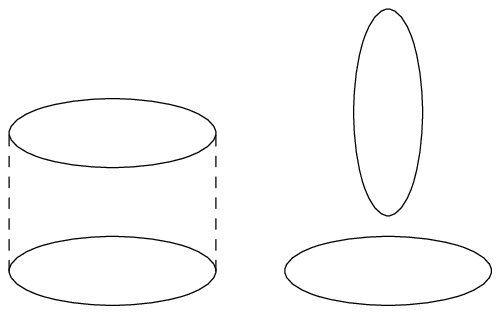}
\caption{Две пары окружностей: изометрически объемлемо однородная и нет}\label{fig:okr-ioo-i-net}
\end{minipage}\hss
\begin{minipage}[b]{\tempa}\centering
\includegraphics{fig4.2.eps}
\caption{Тор вращения в $\R^3$}\label{fig:torR3}
\end{minipage}%
\hss}
\end{figure}

\begin{problem}
%Докажите, что
Следующие подмножества изометрически объемлемо однородны:
\begin{myitemize}
\item пара точек на плоскости;

\item вершины правильного многоугольника на плоскости;

\item целочисленная решетка (т.\,е.~множество точек, все координаты которых целые) на плоскости;
%(и, вообще, в $\R^m$),

\item окружность $S^1:=\{(x,y)\in\R^2\colon x^2+y^2=1\}$ на плоскости;

%\item окружность $S^1:=\{(x,y)\in\R^2\colon x^2+y^2=1\}$ в трехмерном пространстве
%(и, вообще, в $\R^m$ при $m\ge2$),

\item сфера $S^2:=\{(x,y,z)\in\R^3\colon x^2+y^2+z^2=1\}$ в трехмерном пространстве (рис.~\ref{fig:sphere});
%(и, вообще, в $\R^m$ при $m\ge3$),

\item винтовая линия в трехмерном пространстве (рис.~\ref{fig:vint}), т.\,е.~линия, заданная параметрическим уравнением $r(t)=(t,\cos t,\sin t)$;\footnote{По такой кривой движется электрон в постоянном магнитном поле, если напряженность $H$ является постоянным вектором и начальная скорость электрона не параллельна и не перпендикулярна напряженности.
Это можно доказать, используя закон Био"--~Савара"--~Лапласа движения электрона, утверждающий, что
$\ddot\gamma=\dot\gamma\times H$.}

\item объединение двух окружностей в трехмерном пространстве, ограничивающих основания прямого кругового цилиндра (т.\,е.~двух окружностей, одна из которых получена из другой параллельным переносом на вектор, перпендикулярный их плоскостям),
см. рис. \ref{fig:okr-ioo-i-net};

\item тор в $\R^4$ (рис.~\ref{fig:torR4}), являющийся произведением двух окружностей (или, что то же самое, заданный параметрическим уравнением $r(s,t)\hm=(\cos s,\sin s,\cos t,\sin t)$).
\end{myitemize}
\end{problem}

Все эти примеры могут быть тривиально обобщены на высшие
размерности. Действительно, легко сообразить, что если
плоское изометрически объемлемо однородное подмножество
рассмотреть как подмножество трехмерного пространства, то
оно также будет изометрически объемлемо однородным.

%А вот примеры подмножеств, не являющихся изометрически объемлемо однородными.

%\smallskip
%\centerline{\it Рисунок: две пары окружностей: и.о.о. и нет; тор вращения в $\R^3$}
%\smallskip

\begin{problem}
%Докажите, что
Следующие подмножества не являются изометрически объемлемо однородными:
\begin{myitemize}
\item множество вершин неравностороннего треугольника на плоскости;

\item отрезок в $\R^1$ (указание: рассмотрите его крайнюю точку);

\item объединение пересекающихся прямых (указание: рассмотрите точку их пересечения);

\item парабола $y=x^2$ на плоскости;

\item объединение двух окружностей в трехмерном пространстве,
отличное от приведенного в предыдущей задаче (рис. \ref{fig:okr-ioo-i-net});

\item тор вращения в $\R^3$ (рис.~\ref{fig:torR3}; указание:
у школьников может не получиться доказать это).
\end{myitemize}
\end{problem}

Сформулируем естественную гипотезу о характеризации изометрически объемлемо однородных подмножеств.
Для этого нам понадобятся еще два определения.
%Эта часть не используется в дальнейшем и может быть пропущена.

Подмножество плоскости (или трехмерного пространства)
называется \emph{замкнутым}, если для любой точки его
дополнения имеется круг (шар) положительного радиуса с центром в
этой точке, пересечение которого с нашим подмножеством
пусто.

Подмножество плоскости (или трехмерного пространства)
называется \emph{связным}, если на плоскости не существует
двух непересекающихся замкнутых  множеств, пересечение
каждого из которых с нашим подмножеством непусто.

%{\it Непрерывно дифференцируемой кривой} на плоскости называется образ непрерывно
%дифференцируемого отображения $\gamma\colon\R\to\R^2$, для которого скорость $\dot\gamma(t)\ne0$ при любом $t$.
%Аналогично определяются непрерывно дифференцируемые кривые в трехмерном пространстве и в $\R^m$.

\begin{conjecture}\label{isom}
\textup{(a)} Изометрически объемлемо однородное связное
замкнутое подмножество плоскости является точкой, прямой,
ок\-руж\-ностью или всей плоскостью.

\textup{(b)} Изометрически объемлемо однородное связное
замкнутое\linebreak подмножество трехмерного пространства является
точкой, прямой, окружностью, винтовой линией, сферой,
цилиндром или всем пространством.
\end{conjecture}

У  этой гипотезы есть аналог и для $\R^m$.

Эту гипотезу можно легко доказать для подмножеств, являющихся дважды дифференцируемыми кривыми
(их определение аналогично приведенному ниже перед теоремой \ref{thm:t3}) с использованием понятия \emph{кривизны}.
Общий случай можно попытаться доказать с использованием классификации движений (эту идею сообщил мне А.\?Ошемков).
Идея неэлементарного доказательства приведена в \ref{pril}.
%Возможно, эта теорема известна специалистам.

\begin{problem}
Подмножество $N$ плоскости называется \emph{переносно объемлемо однородным}, если для любых двух точек $x,y\in N$ существует параллельный перенос, переводящий $x$ в $y$, а $N$ в себя.

(a) Переносно объемлемо однородное замкнутое связное
подмножество плоскости является точкой, прямой или всей
плоскостью. (Эта задача является шагом к доказательству
вышеприведенной гипотезы, поэтому интересно прямое
доказательство, а не вывод из приведенной ги\-по\-тезы.){\looseness=1\par}

(b)** Верно ли, что переносно объемлемо однородное связное подмножество плоскости является точкой, прямой или всей плоскостью?
%Верна ли гипотеза \ref{isom} без предположения замкнутости?
(Например, может ли какая"=нибудь <<дикая>> подгруппа плоскости по сложению, которая строится с помощью аксиомы выбора, быть связной?)
\end{problem}

\begin{problem} Определите \emph{подобистическую объемлемую однородность} подмножеств плоскости.

(a) Приведите пример подобистически объемлемо однородного
подмножества плоскости, не являющегося изометрически
объемлемо однородным.{\looseness=1\par}

(b)** Попробуйте охарактеризовать связные замкнутые подобистически объемлемо однородные подмножества плоскости.
\end{problem}

\subsection{Аффинная объемлемая однородность}\label{subsec1.2}

%Допустим теперь, что вынимаемый объект не жёсткий, а допускает деформации; скажем, мы хотим
%вытащить из бронированного (и потому жёсткого) кабеля его <<мягкую>> сердцевину.
%Какой формы может быть электрический кабель, чтобы провод можно было вытащить из него?

Какой формы может быть металлический кабель, чтобы из него
можно было вытащить его <<мягкую>> сердцевину? Кабель
деформировать нельзя (он жесткий), а провод можно
деформировать плавно, но нельзя ломать. Математическая
формулировка этого вопроса приводит к понятию \emph{дифференцируемой}
% !!! сделать ссылку
объемлемой однородности из \ref{sec4}
% ??? Зачем вообще приводить название параграфа?
% оставить
<<Обобщение на диффеоморфизмы>>. Мы сначала рассмотрим
более простое понятие \emph{аффинной} объемлемой
однородности. Оно хуже отражает ситуацию, зато доступно
школьнику и интересно с точки зрения математики. А самое
главное, на
\settowidth{\wrfwidth}{\includegraphics[scale=1]{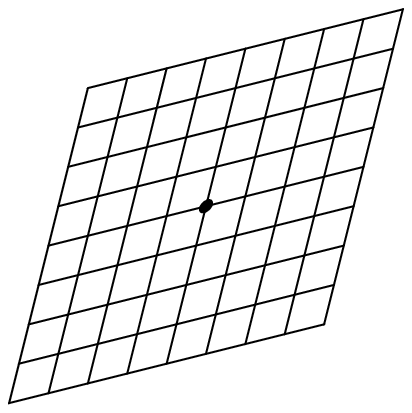}}
\begin{wrapfigure}{o}{0pt}\centering
\includegraphics[scale=1]{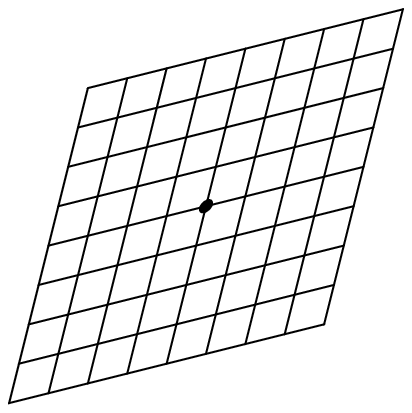}
\caption{Образ решетки при аффинном преобразовании}\label{fig:lataff}
\end{wrapfigure}
примере его изучения в этой брошюре показана
идея доказательства характеризации \textit{дифференцируемо}
объемлемо однородных подмножеств.

%\smallskip
%{\bf Определение.}
%\begin{definition*}
\qdefinitions
Подмножество $N$ пространства $\R^m$ (в частности,
плоскости $\R^2$ или трехмерного пространства $\R^3$)
называется \emph{аффинно объемлемо однородным}, если для
любых двух точек $x,y\in N$ существует
аффинное
преобразование $h\colon\R^m\to\R^m$, переводящее $x$ в $y$, а $N$
в себя.
%\end{definition*}

%\smallskip
%\centerline{\it Рисунок: кошка и ее образ при аффинном преобразовании}

\begin{figure}[t]\centering
\includegraphics{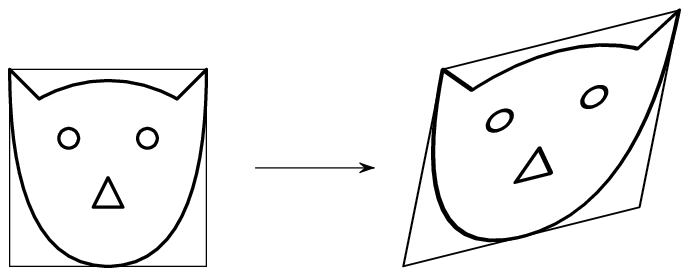}
\caption{\smash[b]{Кошка и ее образ при аффинном преобразовании}}\label{fig:cat-aff}
\end{figure}

%\smallskip
Напомним, что \emph{аффинным преобразованием} плоскости
называется композиция движения, гомотетии и  растяжения
относительно прямой, см. рис. \ref{fig:lattice}, \ref{fig:lataff}
и \ref{fig:cat-aff}. Здесь растяжение относительно прямой
можно заменить на параллельную
проекцию из одной копии
нашей плоскости, находящейся в трехмерном пространстве, на
другую. Подробнее см.~\cite{Pr}.

%\begin{figure}[h]
%\hfill \includegraphics[scale=.5]{transform.3}\quad
%\hfill \mbox{}
%\caption{Образ решетки при аффинном преобразовании}\label{fig:lataff}
%\end{figure}

Аффинное преобразование трехмерного пространства определяется более сложно;
мы приведём здесь это определение, хотя до \ref{sec4} оно нам не понадобится.
Пусть заданы точки $O$ и $O'$, а также две некомпланарные тройки векторов $a,b,c$ и $a',b',c'$.
Тогда \emph{аффинным преобразованием} трехмерного пространства, отвечающим $O,a,b,c$ и $O',a',b',c'$,
называется преобразование, переводящее точку $O+xa+yb+zc$ в точку $O'+xa'+yb'+zc'$.
Аналогично определяется аффинное преобразование пространства $\R^m$; такое определение для $m=2$
равносильно вышеприведенному.

\begin{problem}
Cледующие подмножества плоскости аффинно объемлемо однородны:
\begin{myitemize}
\item любое изометрически объемлемо однородное подмножество;

\item эллипс, заданный уравнением $x^2+2y^2=1$;

\item парабола $y=x^2$;

\item
%(равносторонняя)
гипербола $y=1/x$.
\end{myitemize}
\end{problem}

\begin{figure}[h]\centering
\hbox to\textwidth{\hss
\includegraphics{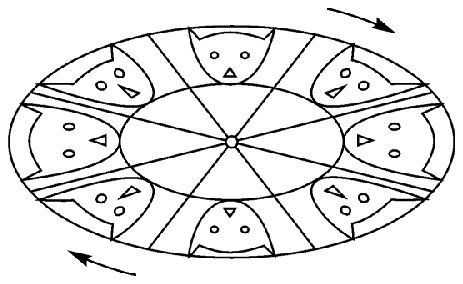}\hss
\includegraphics{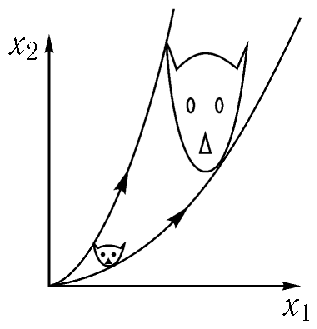}%
\hss}
\includegraphics{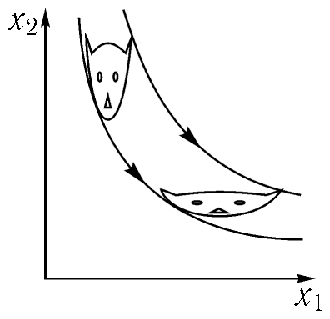}
\caption{Кошка и ее образы при эллиптическом, параболическом и гиперболическом поворотах}\label{cat-ell-par-hyp}
\end{figure}

Подсказка: см. рис.~\ref{cat-ell-par-hyp}.
%\smallskip
%{\it Рисунок: кошка и ее образы при эллиптическом, параболическом и гиперболическом поворотах [Ar84, рис. 132, 46, 47]}

%\smallskip
Какие еще бывают аффинно однородные подмножества плоскости?
%А может ли быть что"=нибудь ещё? К примеру, возьмём график непрерывной функции $f\colon\R\to \R$.

\begin{problem}
(a) Не любое конечное множество на плоскости является аффинно объемлемо однородным.

(b)* Опишите конечные аффинно объемлемо однородные подмножества плоскости.
(Ответ "--- аффинно правильные многоугольники.
Полезно использовать, что аффинное преобразование сохраняет центр масс и <<эллипс инерции>>.)
\end{problem}

\begin{problem}
(a) График функции $y=|x|$ не является аффинно объемлемо однородным подмножеством плоскости.

(b) Если функция $f\colon\R\to \R$ дифференцируема хотя бы в одной точке и график $f$ аффинно объемлемо однороден, то $f$ дифференцируема в любой точке. (Если производная в точке равна плюс бесконечности или равна минус бесконечности, то мы считаем функцию дифференцируемой в этой точке.)
\end{problem}

%Можно показать, что график непрерывной функции $f\colon\R\to\R$, дифференцируемой в одной точке
%и не дифференцируемой в другой, не является аффинно объемлемо однородным.

\begin{figure}[h]\centering
\hbox to\textwidth{\hss
$\vcenter{\hbox{\includegraphics{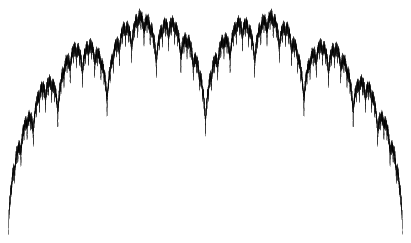}}}$\hss
$\vcenter{\hbox{\includegraphics{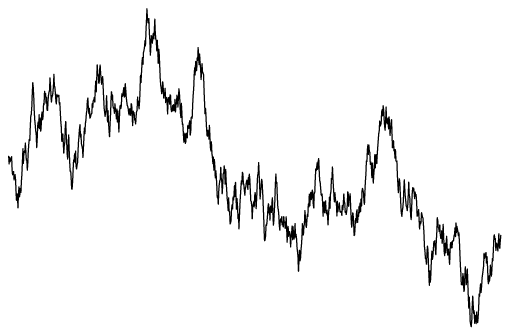}}}$%
\hss}
\caption{Пила Вейерштрасса и броуновское движение}\label{fig:WBC}
\end{figure}

Но ведь есть и непрерывные функции, не дифференцируемые ни в одной точке.
Например, \emph{пила Вейерштрасса} $f(x):=\textsum_{n=0}^\infty2^{-n}\sin(13^n\pi x)$,
см. рис.~\ref{fig:WBC}.
Такие примеры встречаются и в физике при изучении \emph{броуновского движения}.
Что тогда? Может ли непрерывная функция быть <<одинаково не дифференцируемой>> во всех точках, т.\,е.~может ли её график быть аффинно объемлемо однороден?
Оказывается, что нет.

\begin{statement}\label{p:diff-f}
Если график непрерывной функции $f\colon\R\to\R$ аффинно объемлемо однороден, то эта функция дифференцируема в любой точке.
%(Если производная в точке бесконечна, то мы считаем функцию дифференцируемой в этой точке.)
Более того, ее производная непрерывна.\footnote{%
Если производная в точке равна плюс бесконечности или равна
минус бесконечности, то мы считаем функцию дифференцируемой
в этой точке. Непрерывность производной в такой точке
означает, что производная бесконечно большая при стремлении
аргумента к этой точке.}
\end{statement}

%Доказательство этого результата и его обобщений "--- основная цель настоящей брошюры.
Читатель, которому остаток этого пункта и следующий покажутся\linebreak слишком трудными, может сразу перейти к одному из двух следующих параграфов.
%, где доказано утверждение \ref{p:diff-f}.

Утверждение~\ref{p:diff-f} является частным случаем более общего факта, который мы сейчас сформулируем.
Для этого нам понадобятся определение замкнутости (см. выше перед гипотезой \ref{isom}) и следующее определение.
\emph{Непрерывно дифференцируемой кривой} на плоскости называется образ непрерывно  дифференцируемого отображения $\gamma\colon\R\to\R^2$, для которого скорость $\dot\gamma(t)\ne0$ при любом $t$.

\begin{theorem}[\cite{RSS93}]\label{thm:t3}
Аффинно объемлемо однородное замкнутое подмножество плоскости является либо

% ??? заменить точки на буквы-цифры? желательно при этом немного всё перестроить, чтобы
%     убрать эти "либо"
% пока оставим (здесь и далее)
$\bullet$ набором изолированных точек, либо

$\bullet$ объединением изолированных (т.\,е.~имеющих непересекаюшиеся окрестности) непрерывно дифференцируемых кривых, либо

$\bullet$ всей плоскостью.
\end{theorem}

%Более того, то же верно и в старших размерностях (см.???):
%В настоящей брошюре мы разберём ход их доказательства. Впрочем, мы начнём с более слабого %результата, доказательство которого, однако, содержит все необходимые идеи.

\begin{remark}
Подмножество плоскости называется \emph{локально замк\-ну\-тым}, если любая его
точка имеет такую замкнутую окрестность~$U$ в плоскости, что пересечение $U$ с нашим
подмножеством замкнуто.
Условие замкнутости в теореме~\ref{thm:t3} можно ослабить до локальной замкнутости
(ибо в том месте доказательства, где используется замкнутость, достаточно
локальной замкнутости).
Теорема~\ref{thm:t3} неверна без предположения замкнутости (или локальной замкнутости);
контрпримером
%в этом случае
является $\Q\hm\subset\R$.{\looseness=1\par}
\end{remark}

\begin{conjecture}\label{aff}
Любое связное замкнутое
аффинно объемлемо однородное подмножество плоскости является либо точкой, либо прямой,
либо эллипсом, либо ветвью гиперболы, либо параболой, либо всей плоскостью.
(Определение связности приведено выше перед гипотезой \ref{isom}.)
\end{conjecture}

%(d) Предыдущее неверно без предположения связности или замкнутости (даже для множеств, не %являющихся наборами изолированных точек).

%Возможно, утверждение этой гипотезы известно специалистам.

%\smallskip
Как мы уже отметили, утверждение \ref{p:diff-f} вытекает из теоремы \ref{thm:t3}.
При этом теорема \ref{thm:t3} является существенно более сильным результатом, чем утверждение \ref{p:diff-f}. Поясним это примерами.

%http://en.wikipedia.org/wiki/Cantor\_dust,

Когда теорема \ref{thm:t3} еще не была доказана, подмножествами плоскости, подозрительными на аффинную объемлемую однородность, были некоторые <<фракталы>>.
Это подозрение основывалось на том, что в другом, более слабом, смысле, они все"=таки однородны (см. следующий пункт).
Определим эти <<фракталы>> "--- стандартное канторово множество и обобщенное канторово множество.

\begin{figure}[t]\centering
\vskip\topskip\vskip-6pt
\settowidth{\tempa}{\includegraphics{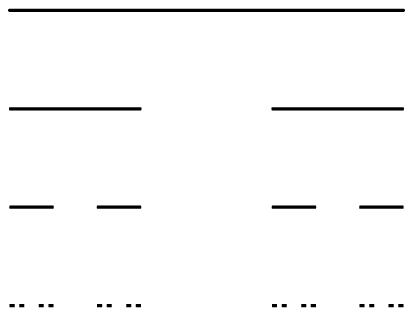}}
\settowidth{\tempb}{\includegraphics{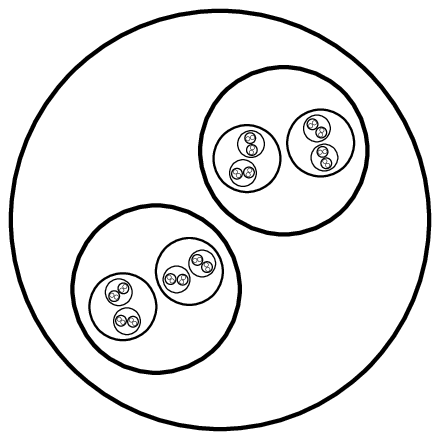}}
\hbox to\textwidth{\hss
\begin{minipage}[b]{\tempa}\centering
\includegraphics{cantor.1.eps}
\caption{Стандартное канторово множество}\label{fig:Cantor-st}
\end{minipage}\hss\hss
\begin{minipage}[b]{\tempb}\centering
\includegraphics{cantor.7.eps}
\caption{Обобщенное канторово множество}\label{fig:Cantor-obob}
\end{minipage}%
\hss}
\vskip12pt
\end{figure}

\emph{Стандартное канторово множество} определяется так:
$$
C:=\biggl\{\textsum_{k=1}^\infty a_k3^{-k}\in[0,1]\colon a_k\in\{0,2\}\biggr\}
$$
(рис.~\ref{fig:Cantor-st}).
%Выше мы уже упоминали о канторовом множестве; напомним, как оно устроено.
Иными словами, стандартное канторово множество
%(именно оно изображено на рис.~\ref{fig:WBC} справа)
получается как предел следующей процедуры: из отрезка $[0,1]$ удаляется средняя треть, затем удаляется средняя треть из каждого из двух полученных отрезков, затем средняя треть из каждого из полученных четырёх отрезков, и так далее.
\begin{figure}[t]\centering
\includegraphics{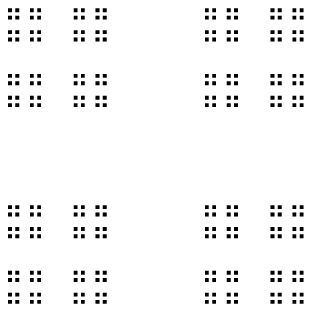}
\caption{Канторова пыль}\label{fig:Cantor-dust}
\end{figure}

%Является ли оно аффинно объемлемо однородным?
%Однако в топологии, говоря, что какое"=либо множество "--- канторово, подразумевают, что оно %устроено так же, как и стандартное канторово, с точки зрения топологии.
Это определение можно обобщить двумя эквивалентными способами (выберите из них наиболе понятный для вас).
Можно определить \emph{обобщенное канторово множество} на плоскости как образ стандартного канторова множества при непрерывном инъективном отображении.
А можно "--- построив <<иерархическую структуру>> следующим образом.

%\smallskip
%{\bf Определение.}
\emph{Обобщенным канторовым множеством} на плоскости называется пересечение объединений
{\advance\abovedisplayskip-3.5pt
\advance\belowdisplayskip-1pt
$$
\bigcap\limits_{n=0}^\infty\,\bigcup\limits_{\alpha_1\dots\alpha_n\in\Z_2^n}
C_{\alpha_1\dots\alpha_n},
$$}%
где $\{C_{\alpha_1\dots\alpha_n}\}$, $n=0,1,2,\dots,$
$\alpha_1\dots\alpha_n\in\Z_2^n$, "--- набор замкнутых непустых подмножеств плоскости,
для которых

$\bullet$ $C_{\alpha_1\dots\alpha_n}
\supset C_{\alpha_1\dots\alpha_n0}\cup C_{\alpha_1\dots\alpha_n1}$
при любых $n$, $\alpha_1\dots\alpha_n$,

$\bullet$ для любого $n$ множества~$C_{\alpha_1\dots\alpha_n}$ попарно не пересекаются и

$\bullet$ $\lim\limits_{n\to\infty}\max\limits_{\alpha_1\dots\alpha_n}
\diam C_{\alpha_1\dots\alpha_n}=0$, где $\diam$ обозначает диаметр множества (рис.~\ref{fig:Cantor-obob}).

%\smallskip
Примером обобщенного канторова множества на плоскости является \emph{канторова пыль} (рис. \ref{fig:Cantor-dust}).
%, ср. \linebreak http://en.wikipedia.org/wiki/Cantor\_dust [засмотрено 20.05.2011]

Из теоремы~\ref{thm:t3} немедленно вытекает

\begin{corollary}\label{cor1}
Никакое обобщенное канторово множество на пло\-с\-кос\-ти не является аффинно объемлемо однородным.
\end{corollary}

В следующем пункте объясняется, почему это следствие выглядит просто чудом
(см. задачи~\ref{prob10}(c,~d)).

%Ковер Серпинского на плоскости не является аффинно объемлемо однородным.

\subsection{Другие виды однородности}\label{s:Cantor}

%Возьмём теперь то определение, с которого мы начали эту брошюру, и оставим в нём только %непрерывность:

%\smallskip
%{\bf Определение.}
\settowidth{\wrfwidth}{\includegraphics[scale=.9]{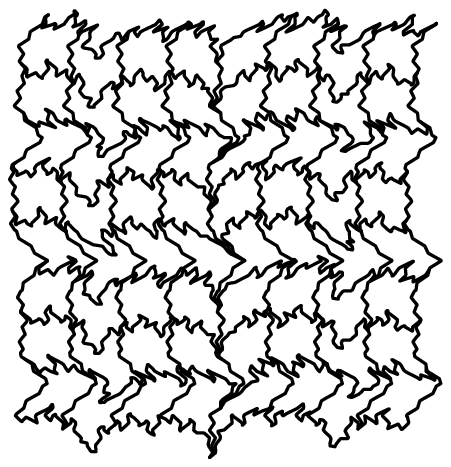}}
\begin{wrapfigure}{o}{\wrfwidth}\centering
\vskip-10pt
\includegraphics[scale=.9]{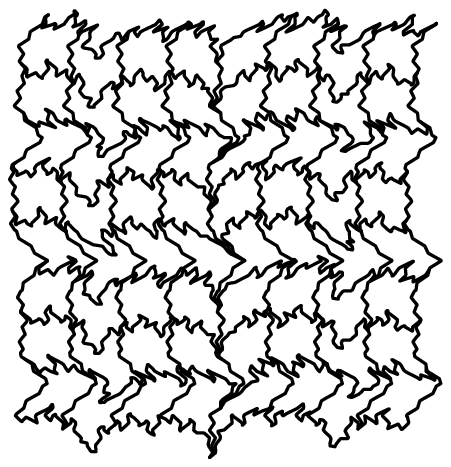}
\caption{Образ решетки при гомеоморфизме}\label{fig:latgomeo}
\end{wrapfigure}
%\begin{definition*}
\qdefinitions
%Замкнутое ограниченное
Подмножество $N$ пространства $\R^m$ (в частности, прямой $\R$ или плоскости $\R^2$)
%или пространства $\R^3$)
называется \emph{однородным}, если
для любых двух точек $x,y\in N$ существует непрерывная биекция
% ??? нужно ли это "на"?
% нужно
(т.\,е.~взаимно однозначное отображение на) $h\colon N\to N$, переводящая $x$ в~$y$,
см. рис. \ref{fig:latgomeo}.
%\end{definition*}

%\begin{figure}[h]
%\hfill \includegraphics[scale=.5]{transform.5.eps}\quad
%\hfill \mbox{}
%\caption{Образ решетки при гомеоморфизме}\label{fig:latdif}
%\end{figure}

%\smallskip
Это определение имеет два существенных отличия от предыдущих.
Во"=первых, отображение $h$ задано только на $N$, а не на объемлющем пространстве $\R^m$
(в отличие от определений \textit{объемлемой} однородности).
Во"=вторых, отображение $h$ предполагается всего лишь непрерывным (а не изометрией, не аффинным~и~т.\,д.).
Поэтому ясно, что любое изометрически или аффинно объемлемо однородное подмножество является однородным.

\begin{figure}[t]\centering
\vskip\topskip\vskip-6pt
\begin{minipage}[b]{.5\textwidth}\centering
\includegraphics[scale=.5]{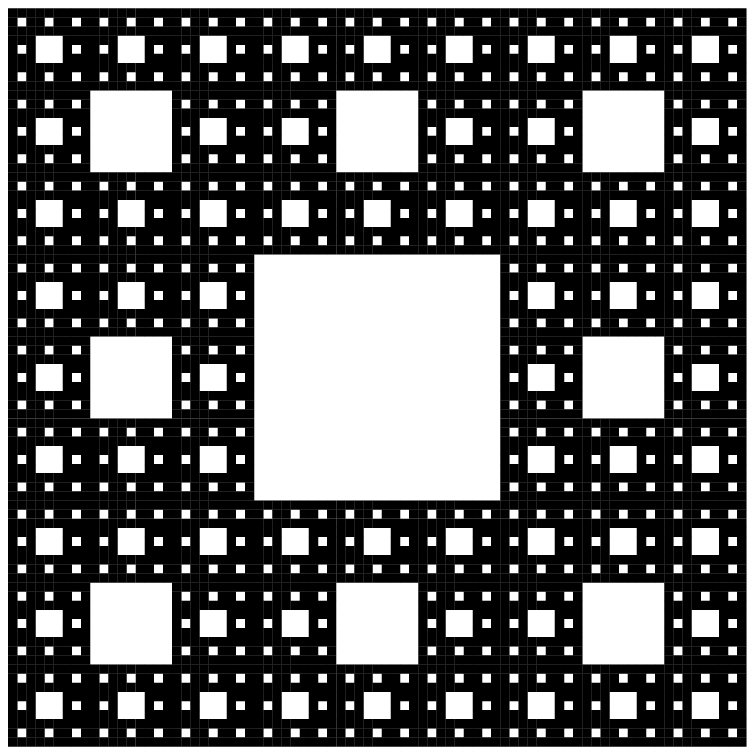}
\caption{Ковер Серпинского}\label{fig:Sierpinski-gasket}
\end{minipage}%
\begin{minipage}[b]{.5\textwidth}\centering
\includegraphics[scale=1.8]{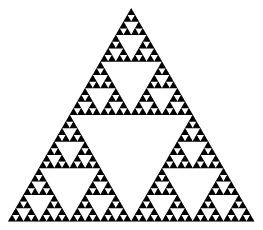}
\caption{Треугольник Серпинского}\label{fig:Sierpinski-triang}
\end{minipage}
\end{figure}

\begin{problem}
Cледующие множества однородны:
\begin{myitemize}
\item конечное множество точек;

\item тор вращения в $\R^3$;

\item множество рациональных точек на отрезке $(0,1)$;

\item стандартное канторово множество\footnote{Некоторые множества из этой задачи определены в предыдущем пункте.};

\item обобщенное канторово множество;

\item \emph{ковер Серпинского} (он строится аналогично канторову множеству, только вместо вырезания средней трети из отрезков происходит вырезание <<сердцевины>> из квадратов, см. рис. \ref{fig:Sierpinski-gasket});

\item орбита непрерывного действия\footnote{А некоторые нет... Напомним, что утверждения, содержащие незнакомые вам термины, можно игнорировать.}
топологической группы на $\R^m$.
\end{myitemize}
\end{problem}

 % http://en.wikipedia.org/wiki/Sierpinski\_carpet

\begin{problem}
Cледующие множества не однородны:
\begin{myitemize}
\item отрезок (указание: рассмотрите его крайнюю точку);
\item объединение пересекающихся прямых (указание: рассмотрите точку их пересечения);
\item \emph{треугольник Серпинского} (рис.~\ref{fig:Sierpinski-triang}; указание: рассмотрите вершины и середины сторон большого треугольника).
\end{myitemize}
\end{problem}

%Для окружности и сферы непрерывная биекция $h\colon N\to N$ из определения однородности может быть %продолжена до аффинного преобразования объемлющего пространства.
%Иными словами, окружность и сфера аффинно объемлемо однородны в $\R^2$ и $\R^3$,
%соответственно.

Следующее понятие формализует свойство жесткого кабеля, необходимое для того, чтобы из него можно было вытащить его мягкую сердцевину, которую можно изгибать и ломать, но нельзя разрывать.

%\smallskip
%{\bf Определение.}
\begin{definition*}
Подмножество $N$ пространства $\R^m$ (в частности, прямой $\R$ или плоскости $\R^2$)
называется \emph{непрерывно объемлемо однородным}, если для
любых двух точек $x,y\in N$ существует непрерывная биекция $h\colon\R^m\to\R^m$,
переводящая $x$ в $y$, а $N$ в себя.
\end{definition*}

\begin{problem}\label{prob10}
(a) Замкнутое ограниченное непрерывно объемлемо однородное\linebreak подмножество прямой
состоит из одной или двух точек.

(b) Канторово множество на прямой, а также ковер Серпинского на плоскости
не являются непрерывно объемлемо однородными.

%(c) Все вышеприведенные примеры однородных множеств, кроме ковра Серпинского, непрерывно объемлемо однородны (в соответствующих евклидовых пространствах).
(c) Однородные множества из задач 8(b,~c,~e,~g) непрерывно объемлемо однородны
(в~соответствующих евклидовых пространствах).

(d) Определите \emph{липшицеву объемлемую однородность} (потребовав от $h$ липшицевости).
Докажите, что обобщенное канторово множество на плоскости является липшицево объемлемо однородным. \end{problem}

%Ввиду задач (c) и (d) следствие \ref{cor1} из предыдущего пункта выглядит просто чудом.

%%%!!!
\newpage
\section{Принцип вложенных отрезков, или\protect\headbreak примени теорему Бэра о категории}\label{s:Baire}

Этот цикл задач посвящен теореме Бэра о категории "--- мощному средству
доказательства теорем существования (подробнее см.~\cite{KF, Ox71}).
В анализе с помощью нее доказывается, например, теорема Банаха об обратном
операторе, которая применяется для доказательства существования решений
нелинейных уравнений.
В топологии теорема Бэра применяется, например, к вложениям компактов и к
аппроксимации отображений гомеоморфизмами.
В этой брошюре мы применим теорему Бэра к доказательству утверждения \ref{p:diff-f} из пункта \ref{subsec1.2} и его обобщений.
%(Формально, при доказательстве утверждения 0 используются только задачи 1 и 3d,
%а при доказательстве обобщений "--- их обобщения с заменой прямой на $\R^m$.)

К задачам приводятся указания в конце этого пункта.
Попытайтесь сначала решить задачи, не читая указаний!

\addvspace{0pt plus1pt}

\begin{problem}\label{b0}
(a) Пусть объединение открытых интервалов $U\subset\R$ неограничено.
Докажите, что существует такое $x$, что
$nx\in U$ для бесконечно большого количества целых $n$.

(b) Дана бесконечно дифференцируемая функция $f\colon\R\to\R$,
%причем $f^{(n)}(x)=0$ для любого $x$ и для всех чисел $n>N_x$.
причем для любого $x$ существует такое целое $N_x$, что $f^{(n)}(x)=0$ для любого $n>N_x$.
Докажите, что $f$ "--- многочлен.
\end{problem}

\addvspace{0pt plus1pt}

Напомним, что подмножество $U\subset\R$ называется\nopagebreak

$\bullet$ \emph{открытым}, если для любого $x\in U$ существует такое $\varepsilon>0$, что
$(x-\varepsilon,x+\varepsilon)\subset U$;

$\bullet$ \emph{всюду плотным}, если для любых $a,b\in\R$ пересечение $(a,b)\cap U$ непусто.

\addvspace{0pt plus1pt}

\begin{problem}\label{b1}
\textit{Теорема Бэра о категории}. Докажите, что пересечение счетного числа открытых всюду
плотных подмножеств прямой является всюду плотным (и, в частности, непустым).
\end{problem}

\addvspace{0pt plus1pt}

%Докажите, что существует такая последовательность $x_n\in\R$, что для любой
%возрастающей ограниченной последовательности $y_n\in\R$ найдется номер $n$,
%для которого $|x_n-y_n|\le1/n$.
%Без т. Бэра: берем $|x_n-x_{n+1}|=\frac1n$,
%$x_n$ не ограничена ни сверху, ни снизу

\begin{problem}\label{b2}
Докажите, что если функция двух переменных непрерывна по каждой
переменной, то она имеет точку непрерывности.
\end{problem}

\addvspace{0pt plus1pt}

\begin{problem}\label{b3}
Докажите следующее для функций $\R\to\R$.

(a) Поточечный предел последовательности $f_n$ непрерывных функций
(т.\,е.~функция $f(x):=\lim\limits_{n\to\infty}f_n(x)$) обязательно имеет точку непрерывности.

(b) Производная любой дифференцируемой функции имеет точку непрерывности.
\end{problem}

\addvspace{0pt plus1pt}

\begin{problem}\label{b4}
 Для бесконечно дифференцируемой функции $f\colon\R\to\R$, любого $x$ и
%\newline
%
(a) бесконечной последовательности чисел $n$ (зависящей от $x$);
(b) некоторого $n=n_x$
%
%\noindent
выполнено $f^{(n)}(x)=0$. Докажите, что $f$ "--- многочлен.
\end{problem}

\addvspace{0pt plus1pt}

\begin{problem}\label{b5}
(a) Докажите, что прямая не представима в виде объединения попарно
непересекающихся замкнутых отрезков, каждый из которых отличен от точки.

(b) Докажите, что плоскость не представима в виде объединения замкнутых кругов с попарно
непересекающимися непустыми внут\-рен\-нос\-тями.
\end{problem}

\addvspace{0pt plus1pt}

 \begin{problem}\label{b6}
(a) Дано замкнутое ограниченное подмножество $A\subset\R^2$.
Известно, что для любых двух точек $x,y\in A$ существует разбиение $A=X\sqcup Y$ на замкнутые
множества, для которого $x\in X$ и $y\in Y$ (такие множества называются
\emph{нульмерными}).
Докажите, что существует непрерывное инъективное отображение
(т.\,е.~\emph{вложение} или реализация) $a\colon A\to\R$.

(b)* Дано замкнутое ограниченное подмножество $A\subset\R^{100}$.
Известно, что для любых двух точек $x,y\in A$ существует разложение $A=X\cup Y$ в объединение
замкнутых множеств, пересечение которых нульмерно, причем $x\in X$ и $y\in Y$
(такие множества называются \emph{одномерными}).
Докажите, что существует непрерывное инъективное отображение
%(т.\,е.~вложение или реализация)
$a\colon A\hm\to\R^3$.{\looseness=1\par}

(c)* \textit{Теорема Менгера"--~Небелинга"--~Понтрягина}.
Дайте определение $n$"=мерного (замкнутого ограниченного) множества в $\R^N$ и
докажите, что любое $n$"=мерное множество вложимо в $\R^{2n+1}$.
\end{problem}

\subsection*{Указания}
%\smallskip
%\newpage
%\small
%{\bf }

\textbf{\ref{b0}(a).} Сначала докажите, что
\textit{существует такое $x_1\in(0,1)$, что
$n_1x_1\in U$ для некоторого $n_1>1$}.

Тогда \textit{существует такое $\epsilon_1>0$,  что
$n_1(x_1-\epsilon_1,x_1+\epsilon_1)\subset U$}.

Потом докажите, что
\textit{существует такое $x_2\in(x_1-\epsilon_1,x_1+\epsilon_1)$,
что $n_2x_2\in U$ для некоторого $n_2>2$}.

И~т.\,д.

%\smallskip
Такие решения, основанные на принципе вложенных отрезков, удобно придумывать
и записывать на языке теоремы Бэра о категории.
Вышеприведенное решение коротко записывается так:
по теореме Бэра о категории
$\bigcap\limits_{n=1}^\infty\bigcup\limits_{k=n}^\infty\frac1kU\ne\emptyset$.

%\smallskip
\textbf{\ref{b0}(b).} Сначала докажите следующий факт.
{\itshape Пусть $U_n\subset\R$ "--- непустые открытые множества,
$U_1\supset U_2\supset\dots$ и
$\bigcap\limits_{n=1}^\infty U_k=\emptyset$.
Тогда существуют такие $n$ и
%максимальный
интервал $(a,b)\subset U_n$, что

$\bullet$ $(a,b)$ максимален, т.\,е.~$U_n$ не содержит никакого большего
интервала $(c,d)\supset(a,b)$, и

$\bullet$ один из интервалов $(a,a+\epsilon)$ и $(b-\epsilon,b)$
не пересекает множество $U_{n+1}$ для некоторого $\epsilon>0$. }

\textit{Указание к доказательству факта.}
Предположим противное.
Тогда для любого $n$, любого максимального интервала $(a,b)\subset U_n$
и любого $\varepsilon>0$ оба интервала $(a,a+\varepsilon)$ и
$(b-\varepsilon,b)$ пересекаются с $U_{n+1}$.
Теперь докажите, что для любого максимального интервала $(a,b)\subset U_n$
либо $U_{n+1}$, либо $U_{n+2}$ содержит интервал, замыкание которого лежит в
$(a,b)$.
Выведите отсюда, что $\bigcap\limits_{n=1}^\infty U_k\ne\emptyset$.
Противоречие.

\textit{Окончание решения задачи \textup{\ref{b0}(b)}.}
Пусть $f$ не многочлен.
Положим
$U_n:=\R-\bigcap\limits_{k=n}^\infty(f^{(k)})^{-1}(0)$.
Применим приведенный факт.
Получим такие $n$ и максимальный интервал $(a,b)\subset U_n$, что для
некоторого $\epsilon>0$ (не уменьшая общности)
$(a,a+\varepsilon)\cap U_{n+1}=\emptyset$.
Тогда $f^{(n)}(a)=0$ и $f^{(n+1)}(a,a+\varepsilon)=0$.
Поэтому $f^{(n)}[a,a+\varepsilon)=0$. Противоречие.

%\smallskip
\textbf{\ref{b1}.} Это простое следствие принципа вложенных отрезков.

%\smallskip
\textbf{\ref{b2}.}
Множество точек непрерывности функции $f$ "--- это
%\linebreak
$$\bigcap\limits_{n=1}^\infty\bigcup\limits_{k=1}^\infty
\Bigl\{x\colon |fy_1-fy_2|<\frac1n\mbox{ при }x-\frac1k<y_1<y_2<x+\frac1k\Bigr\}.$$

%\smallskip
\textbf{\ref{b3}.} (a)
Фиксируем $\varepsilon>0$.
Положим $U_n:=\bigcup\limits_{i,j\ge n}\{x\colon |f_ix-f_jx|>\varepsilon\}$.
Тогда $U_n$ открыто, $U_n\supset U_{n+1}$ и
$\smash[t]{\bigcap\limits_{n=1}^\infty U_n}=\emptyset$.
Значит, по теореме Бэра для любого отрезка $[a,b]$ существуют число $n$ и интервал
$(c,d)\subset[a,b]$, не пересекающийся с $U_n$.
Значит, $|f_ix-f_jx|\le\varepsilon$ для любых $x\in(c,d)$ и $i,j\ge n$.
Поэтому $|fx-f_nx|\le\varepsilon$ для любого $x\in(c,d)$.

(b) Используйте (a) и $f'(x)=\lim\limits_{n\to\infty}n\Bigl(f\Bigl(x+\frac1n\Bigr)-f(x)\Bigr)$.

\section{Доказательство теоремы~\ref{thm:t3} и утверждения~\ref{p:diff-f}}\label{s:proof-aff}

Здесь мы докажем утверждение \ref{p:diff-f}.
Доказательство теоремы~\ref{thm:t3} аналогично.
При доказательстве можно вместо определения
аффинного преобразования использовать только следующие его свойства.

\textit{Для любого аффинного преобразования $h\colon\R^2\to\R^2$ образ любого
открытого треугольника с вершиной в любой точке $x$ содержит
некоторый открытый треугольник\footnote{Здесь <<содержит некоторый открытый треугольник>> можно было бы заменить на <<является открытым
треугольником>>. Однако для дальнейших обобщений нам удобно сформулировать это свойство именно в приведенном виде.}
с вершиной в точке $h(x)$.}

\textit{Если пересечение аффинно объемлемо однородного подмножества
плоскости с некоторым
%открытым
кругом является непрерывной кривой, имеющей
точку дифференцируемости, то это подмножество является объединением
изолированных дифференцируемых кривых.}

%\smallskip
%\textit{Доказательство утверждения \ref{p:diff-f}.}
\begin{proof}[Доказательство утверждения \ref{p:diff-f}]
Обозначим через $N$ график данной функции $f$.
Возьмем точку $a\in\R^2-N$.
Расстояние от $a$ до $N$ не равно нулю.
Значит, существует точка $y\in N$, для которой $|a-y|$ равно этому расстоянию.
Тогда открытый круг $D$ с центром в $y$ радиуса $|a-y|$ не пересекает $N$.

Обозначим через $R^\varphi$ поворот плоскости на угол
$\varphi$ вокруг начала координат. Обозначим через $B_l$
равнобедренный треугольник (двумерный открытый) с вершиной
в начале координат, углом $2\pi/l$ при вершине и высотой
длины $1/l$, параллельной оси $Oy$.

При любом $x\in N$ существует аффинное преобразование $h\colon\R^2\to\R^2$,
переводящее $y$ в $x$, а $N$ в себя.
Так как $h$ аффинно, то $h(D)\supset x+R^\varphi B_l$ для некоторых
$l$ и $\varphi$.
Поэтому
$$
\advance\thickmuskip-1mu\advance\medmuskip-1mu
\quad \text{при любом $x\in N$ существуют такие $l$ и $\varphi$, что }
(x+R^\varphi B_l)\cap N=\emptyset.\hskip0pt\leqno (*)
$$

%$(*)$\quad \textit{при любом $x\in N$ существуют такие $l$ и $\varphi$, что
%$(x+R^\varphi B_l)\cap N=\emptyset$.}

Возьмем произвольную последовательность $\{\varphi_l\}$, всюду плотную на
$[0,2\pi]$.
Обозначим
$$N_l:=\{x\in N\colon (x+R^{\varphi_l}B_l)\cap N=\emptyset\}.$$
Ввиду условия~$(*)$ имеем $N=\bigcup\limits_{l=1}^\infty N_l$.

% ??? это уже было (слово в слово) тремя страницами ранее
% так и надо
Напомним следующие определение и теорему. Подмножество $U\subset\R$ называется

$\bullet$ \emph{открытым}, если для любого $x\in U$ существует такое $\varepsilon>0$, что
$(x-\varepsilon,x+\varepsilon)\subset U$;

$\bullet$ \emph{всюду плотным}, если для любых $a,b\in\R$ пересечение
$(a,b)\cap U$ непусто.

%\smallskip
%{\bf Теорема Бэра о категории.} {\it
\begin{theoremB}
Пересечение счетного числа открытых всюду
плотных подмножеств прямой является всюду плотным (и, в частности, непустым).
\end{theoremB}
%}

%\smallskip
Нетрудно проверить, что \textit{$N_l$ замкнуто в $N$} (докажите или найдите детали в
\cite[лемма 3.1]{RSS96}).
Значит, по теореме Бэра о категории
%(см. следующий пункт)
некоторое $N_l$
содержит непустое открытое в $N$ множество.
\begin{figure}[t]\centering
\includegraphics{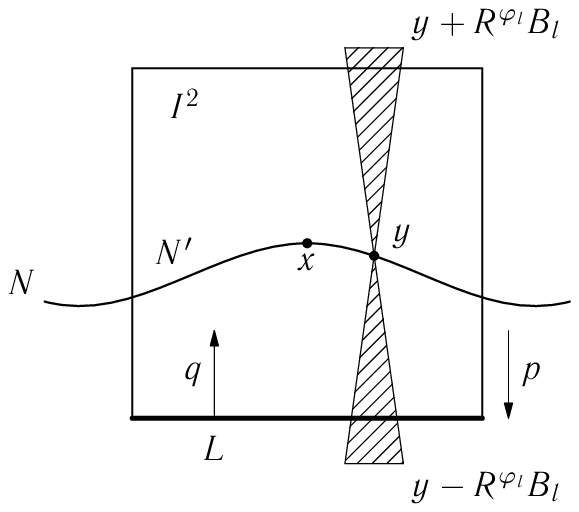}
% ??? сделать подпись к рисунку?
% не надо
\caption{}\label{ris.os}
\end{figure}
%\centerline{\it Рисунок } %approximately here}
%\smallskip

Поэтому существуют точка $x\in N$ и замкнутый квадрат $I^2$ со стороной меньше
$1/l$ с центром в $x$, для которых $N':=N\cap I^2\subset N_l$
(см. рис. \ref{ris.os}%, $A_l:=R^{\varphi_l}$
).
Тогда
$$[(y+R^{\varphi_l}B_l)\cup(y-R^{\varphi_l}B_l)]\cap N'=\emptyset
\quad\mbox{при любом } y\in N'.
\leqno (**)
$$
Действительно, если $z\in(y-R^{\varphi_l}B_l)\cap N'$, то
$y\in(z+R^{\varphi_l}B_l)\cap N'\subset N_l$, что невозможно.

%\begin{figure}
%\hfill        \includegraphics{os-1.pdf} \hfill \mbox{}
%\caption{}\label{fig:Cantor}
%\end{figure}

Можно считать, что угол между некоторой стороной $L$ квадрата $I^2$
и осью $Ox$ равен $\varphi_l$.
Можно также считать, что $N'$ связно и гомеоморфно отрезку (иначе можно заменить $N'$
на малую окрестность точки $a\in N'$, которая гомеоморфна отрезку,
поскольку $N$ "--- график функции).
Тогда ортогональная проекция множества $N'$ на $L$ содержит некоторый отрезок
ненулевой длины.
Можно считать, что этот отрезок совпадает с $L$ (иначе уменьшим $L$).

Напомним, что отображение $q\colon L\to[0,1]$ называется \emph{липшицевым}, если
существует такое $s$, что $|q(x)-q(y)|<s|x-y|$ для любых двух различных точек
$x,y\in L$.
Из~$(**)$ следует, что $N'$ есть график некоторой липшицевой функции
$q\colon L\to[0,1]$ (при естественном представлении $I^2=L\times[0,1]$).
Функция $q$ имеет точку дифференцируемости\footnote{Это следует из того, что любая липшицева функция раскладывается в разность монотонных и
что любая монотонная функция имеет точку дифференцируемости.
Первое несложно, а второе доказывается с использованием соображений меры.
Детали нетривиальны и приведены, например, в~\cite{KF}.
}.
Значит, и данная функция $f$ имеет точку дифференцируемости.
Тогда из аффинной объемлемой однородности вытекает, что $f$
дифференцируема в любой \mbox{точке}.

Докажем теперь часть <<более того>> утверждения \ref{p:diff-f}.
Производная $f'$ имеет точку непрерывности (см. предыдущий параграф).
Тогда из аффинной объемлемой однородности вытекает, что $f$
непрерывно дифференцируема.
%QED
\end{proof}

\section{Обобщение на диффеоморфизмы}\label{sec4}

Напомним, что в \ref{sec1} мы рассматривали жесткий кабель,  из которого можно вытащить его мягкую сердцевину, которую можно изгибать, но нельзя ломать.
Мы обещали ввести понятие, которое более адекватно формализует необходимое свойство такого кабеля.

%\smallskip
%{\bf Определение.}
\begin{definition*}[\cite{DRS89}]
Подмножество $N$ пространства $\R^m$ (в частности,
плоскости $\R^2$ или трехмерного пространства $\R^3$)
называется \emph{дифференцируемо объемлемо однородным}, если
для любых двух точек $\hsk{3}x,y\hm\in N$ существует диффеоморфизм
пространства $\R^m$, переводящий $x$ в $y$, а $N$ в себя.
\end{definition*}

%\smallskip

\settowidth{\wrfwidth}{\includegraphics[scale=.8]{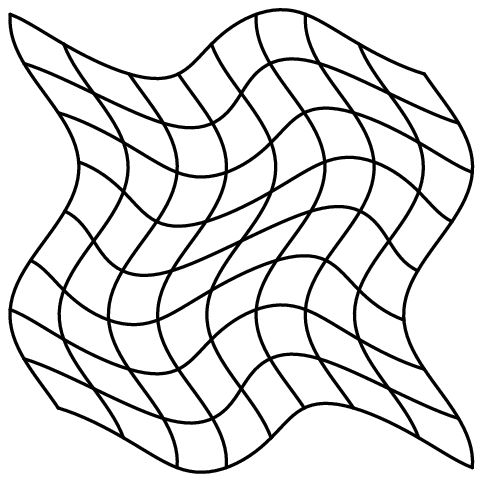}}
\begin{wrapfigure}[14]{o}{\wrfwidth}\centering
\includegraphics[scale=.8]{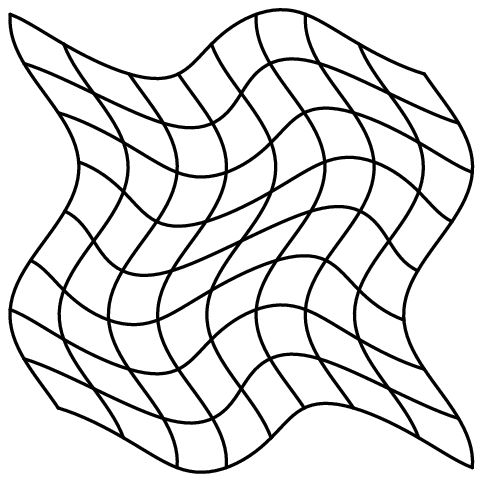}
\caption{Образ решетки при диффеоморфизме}\label{fig:latdif}
\end{wrapfigure}
Диффеоморфизм "--- это локально"=аффинное преобразование.
Формально, взаимно однозначное соответствие $F\colon\R^m\to\R^m$ называется
\emph{диффеоморфизмом}, если для любой точки $z_0\in\R^m$ существует такое
%невырожденное линейное отображение
аффинное преобразование $A\colon\R^m\to\R^m$, что $A(z_0)=F(z_0)$
и отображение $\alpha(z)\hm{:=}\frac{F(z)-A(z)}{|z-z_0|}$ бесконечно малое при
$z\to z_0$, см. рис. \ref{fig:latdif}.
% и бесконечно малая функция $\alpha\colon\R^2\to\R^2$,  и для любой
%точки $z\in\R^2$ выполнено

%\begin{figure}[h]\centering
%\includegraphics[scale=.5]{transform.4.eps}
%\caption{Образ решетки при диффеоморфизме}\label{fig:latdif}
%\end{figure}

\textit{Примерами} диффеоморфизмов плоскости являются большинство геометрических
преобразований: параллельный перенос, поворот, осевая симметрия, гомотетия,
растяжение от прямой.
(Инверсия и центральная проекция являются диффеоморфизмами плоскости без точки
и плоскости без прямой, соответственно.)
Поэтому аффинно объемлемо однородное подмножество является дифференцируемо объемлемо однородным.

Диффеоморфизмом плоскости \textit{не является}, например, отображение, заданное формулой $(x,y)\mapsto(x^3,y)$.
(В самом деле, <<линейная часть>> при $x=0$ не является аффинным преобразованием, ибо не является обратимым отображением.)

\begin{problem} Подмножества плоскости из теоремы \ref{thm:t3} являются дифференцируемо объемлемо однородными.
\end{problem}

Естественно возникает следующий вопрос: какие еще бывают дифференцируемо объемлемо однородные (допустим, замкнутые) подмножества  плоскости?

%Понятно, что нам подходит конечный набор точек и замкнутая гладкая несамопересекающаяся кривая %(или набор таких кривых, ограничивающих попарно непересекающие области).
%Из теоремы \ref{thm:t3} и предыдущей задачи не??? вытекает следующий результат
%Вообще говоря, строгое определение диффеоморфизма опирается на понятия дифференцируемости функций %нескольких переменных и потому выходит за рамки школьной программы. Однако мы надеемся, что оно в %достаточной мере интуитивно понятно. (...)

\begin{theorem}[\cite{RSS96}]\label{dif2}
%Дифференцируемо объемлемо однородное локально замкнутое подмножество плоскости
%является либо набором изолированных точек, либо объединением изолированных дифференцируемых кривых, либо всей плоскостью.
Дифференцируемо объемлемо однородное замкнутое подмножество плоскости
является либо

$\bullet$ набором изолированных точек, либо

$\bullet$ объединением изолированных непрерывно дифференцируемых кривых,
либо

$\bullet$ всей плоскостью.
\end{theorem}

Доказательство аналогично доказательству теоремы \ref{thm:t3}, поскольку свойства,
сформулированные в начале \ref{s:proof-aff}, справедливы для диффеоморфизмов.

По теореме~\ref{dif2} аналоги утверждения \ref{p:diff-f} и следствия \ref{cor1} справедливы с заменой
\textit{аффинной} объемлемой однородности на \textit{дифференцируемую}.

\subsection*{Следствие о диффеоморфизмах прямой}

Взаимно одно\-знач\-ное соответствие $F\colon\R\to\R$ называется
\emph{диффеоморфизмом}, если и оно, и его обратное имеют конечную ненулевую производную в каждой точке.
Примерами диффеоморфизмов прямой являются линейные преобразования: параллельный перенос, растяжение.
Диффеоморфизмом прямой \textit{не является}, например, отображение, заданное формулой
$x\mapsto x^3$.

\begin{corollary}\label{cordif}
Если имеется семейство диффеоморфизмов $h_t\colon\R\hm\to\R$, $t\in\R$,
непрерывно зависящих от параметра $t$, причем $h_s\circ h_t\hm=h_{s+t}$
и $h_{-t}=(h_t)^{-1}$ при любых $s,t\in\R$,  то $h_t$ дифференцируемо
зависит от параметра $t$, т.\,е.~$h_t(x)$ дифференцируемо по $t$
для любого фиксированного $x\in\R$.
\end{corollary}

Обобщение и доказательство см. в следствии \ref{corlie}, приведенном в следующем пункте.

Закончим этот пункт научной переформулировкой и историей следствия~\ref{cordif}, не обязательными для понимания дальнейшего.

Cемейство диффеоморфизмов $h_t\colon\R\to\R$, $t\in\R$, для которых $h_s\circ h_t\hm=h_{s+t}$ и $h_{-t}=(h_t)^{-1}$ при любых $s,t\in\R$, называется \emph{действием} группы $\R$ на прямой (т.\,е.~на себе) диффеоморфизмами.
Такое действие является математическим эквивалентом физического понятия <<двусторонне детерминированный процесс>>; обсуждение  этого можно найти в~\cite[4.2]{Ar84}.
Научная формулировка следствия \ref{cordif}: непрерывное действие группы $\R$ на прямой
диффеоморфизмами является гладким~\cite[4.3]{Ar84}.
%т.\,е.~однопараметрическая группа диффеоморфизмов прямой, net!!!

Около 1986~г. И.\?В.\?Ященко после лекции В.\?И.\?Арнольда спросил у лектора, как доказывать это утверждение. По словам Ященко, Арнольд не смог сходу придумать доказательство и объявил это утверждение открытой проблемой. В процессе обсуждения Арнольдом и Е.\?В.\?Щепиным было открыто (или переоткрыто) понятие дифференцируемой объемлемой однородности, что и привело к доказательству приводимых результатов.

\subsection*{Замечание о других видах однородности}

Напомним, что отображение $q\colon\R^k\to\R^m$ называется \emph{липшицевым}, если
существует такое $s$, что $|q(x)-q(y)|<s|x-y|$ для любых двух различных точек
$x,y\in\R^k$.
Канторово множество может быть \textit{липшицево} объемлемо
однородно вложено в плоскость (докажите или см.~\cite{MR99}).
Значит, \textit{неверен аналог теоремы \ref{dif2} для липшицевой категории}
(т.\,е.~аналог, полученный заменой дифференцируемой объемлемой
однородности на липшицеву и дифференцируемых кривых на липшицевы).
Из этого же примера вытекает, что
\textit{неверен аналог теоремы \ref{dif2} для непрерывной категории}.{\looseness=1\par}

Напомним, что функция называется $C^r$"=дифференцируемой, если ее $r$"=я
производная существует и непрерывна.
Оказывается, аналог теоремы \ref{dif2} для $C^r$"=категории верен.

\begin{theorem}
Аналог теоремы \ref{dif2} верен для $C^r$"=категории при $r\ge1$.
\end{theorem}

Для $r=1$ редукция этого результата к теореме \ref{dif2} аналогична окончанию
доказательства утверждения \ref{p:diff-f}.
Для $r\ge2$ доказательство более сложно \cite[theorem B]{Wi08}.

\begin{conjecture}
Аналог теоремы \ref{dif2} верен для аналитической категории.
\end{conjecture}

\iffalse
Напомним, что функция называется $C^r$"=дифференцируемой, если ее $r$"=я
производная существует и непрерывна.
Оказывается, аналог теоремы \ref{dif2} для $C^1$"=категории верен.
Его доказательство аналогично.
%\footnote{В самом конце доказательства надо дополнительно воспользоваться
%тем, что у производной дифференцируемого отображения есть точка непрерывности.
%Тогда из объемлемой $C^1$"=однородности будет вытекать, что $N'$ есть
%$C^1$"=кривая.}

%\smallskip
%{\bf Гипотеза 3.} {\it
\begin{conjecture}
Аналог теоремы \ref{dif2} верен для $C^r$"=категории при
$r\ge2$ (и для аналитической категории).
\end{conjecture}
%}

%\smallskip
Этим вопросом занималась, в частности, А.\?Вилкинсон~\cite{Wi}.
\fi

%\footnote{
%Вопреки~\cite{RSS96}, я не знаю доказательства гипотезы.
%Для этого доказательства можно вместо конусов $B^{m,k}_l$ пытаться
%рассматривать объекты, более гладко втыкающиеся в начало координат, а вместо
%$(m-k)$"=мерных шаров $D\subset L$ (из второго случая) "--- фигуры, имеющие
%больший порядок касания с $p(N')$.
%См. также~\cite{Wi}.

\section{Приложение: обобщение на многомерный\protect\headbreak случай и многообразия}\label{pril}

Понятия и результаты, с которыми мы работали, допускают естественное обобщение на многомерный случай и на случай многообразий.

%\smallskip
%{\bf Определение.}
\begin{definition*}
Подмножество $N$ дифференцируемого мно\-го\-об\-ра\-зия $M$ называется
\emph{дифференцируемо объемлемо однородным}, если для любых двух точек
$x,y\in N$ существует диффеоморфизм многообразия $M$, переводящий $x$ в $y$.
(Не предполагается ни непрерывности производной диффеоморфизма $h$, ни непрерывной
зависимости $h$ от $x,y$.)
\end{definition*}

%\smallskip
Напомним, что подмножество $N\subset M$  дифференцируемого многообразия $M$ называется \emph{дифференцируемым подмно\-го\-об\-ра\-зием},
если для любой точки $x\in N$ найдутся ее окрестность $U$ и диффеоморфизм $U\to\R^k\hm\times \R^{m-k}$, под действием которого пересечение $U\cap N$ переходит в график некоторой дифференцируемой функции $q\colon\R^k\to \R^{m-k}$.
(Это определение, удобное для доказательства нижеследующей теоремы \ref{thm:t5},
равносильно стандартному~\cite{Pr04}.)

Например, дифференцируемо объемлемо однородным является любое дифференцируемое
подмногообразие дифференцируемого многообразия.\linebreak
% (в частности, \textit{график любой дифференцируемой функции из $\R$ в $\R$}).
Следующая теорема показывает, что верно и обратное.
%В этой заметке мы приводим доказательство обратного утверждения.

\begin{theorem}[\cite{RSS96}]\label{thm:t5}
Дифференцируемо объемлемо однородное замкнутое подмножество дифференцируемого
мно\-го\-об\-ра\-зия является дифференцируемым подмно\-го\-об\-ра\-зием.
\end{theorem}

Известно, что многообразия однородны и что однородное пространство
не обязано быть многообразием (пример: канторово множество).
Теорема \ref{thm:t5} показывает, что свойство быть \textit{дифференцируемым}
подмногообразием равносильно \textit{дифференцируемой объемлемой} однородности.
Ср.~\cite{Gl68}.{\looseness=1\par}

%Если в каком"=то из приведенных ниже замечаний или следствий встретятся
%непонятные читателю термины, то его можно опустить без ущерба для понимания
%остального материала.

\subsection*{Следствия о группах Ли}

При помощи теоремы \ref{thm:t5} удобно доказывать, что некоторые группы
являются группами Ли.
Например, из нее вытекает теорема Картана о том, что \textit{любая замкнутая
подгруппа группы Ли является подгруппой Ли}.

\begin{proof}[Идея доказательства гипотез \ref{isom} и \ref{aff}]
(Эта идея сообщена С.\?В.\?Ивановым.)
Группа самосовмещений нашего множества является замкнутой подгруппой группы
движений (или аффинных преобразований), а значит, подгруппой Ли.
Так как множество связно, то ее компонента единицы тоже действует на множестве
транзитивно.
Связные подгруппы Ли соответствуют подалгебрам Ли.
Для данных конкретных случаев подалгебры
%(с точностью до сопряжения - нет!)
можно перечислить.
Орбиты действий соответствующих подгрупп и будут всеми однородными
множествами. Они автоматически будут гладкими подмногообразиями, поэтому
их можно перечислять и с помощью дифференциальной геометрии.
\end{proof}

%Приведенная теорема позволяет свести следующий результат
%\cite[Theorem 3 on p. 208"--~209]{MZ55} к его простому случаю $m=1$ (т.\,е.
%к уравнению Коши $h(s+t)=h(s)+h(t)$):

\begin{corollary}\label{corlie}
Непрерывное действие группы Ли диффеоморфизмами на дифференцируемом мно\-го\-об\-ра\-зии является дифференцируемым,
ср. \textup{\cite[Theorem 3, p. 208"--~209]{MZ55}}.
\end{corollary}

%\textit{Любая однопараметрическая группа $\{h^t\}_{t\in\R}$ диффеоморфизмов
% мно\-го\-об\-ра\-зия, непрерывно зависящих от параметра
%$t$, на самом деле дифференцируемо зависит от этого параметра}.
%теорема позволяет свести этот результат к случаю $m=1$.

%{\it Доказательство.}
\begin{proof}
Обозначим действие через $h^t\colon M\to M$, $t\in G$.
Определим отображение $\gamma\colon G\to G\times M$ формулой $\gamma(t)=(t,h^t(*))$,
где $*\in M$ "--- произвольная точка.
Для любого $a\in G$ отображение $(x,y)\mapsto(xa,h^a(y))$ определяет диффеоморфизм многообразия $G\times M$.
%плоскости.
Он переводит $\gamma(G)$ в себя.
Значит, по теореме \ref{thm:t5} $\gamma(G)$ "--- дифференцируемое подмногообразие.
Поэтому $h^t(*)$ дифференцируемо по $t$.
%Итак, $h^t(*)$ дифференцируемо по $t$ для любого фиксированного $*\in M$.
%QED
\end{proof}

\begin{corollary}\label{hsc}
Если локально компактная топологическая группа эффективно действует на
гладком мно\-го\-об\-ра\-зии диффеоморфизмами, то это группа Ли.
\end{corollary}

%{\it Доказательство.}
\begin{proof}
Ясно, что любая орбита некоторого непрерывного действия топологической
группы на гладком мно\-го\-об\-ра\-зии диффеоморфизмами является гладко
объемлемо однородной.
Топологическая группа $A_p$ $p$"=адических чисел гомеоморфна канторову
множеству.
Поэтому из теоремы \ref{thm:t5} вытекает, что

\textit{$A_p$ не может свободно (и даже эффективно) непрерывно действовать на
гладком многообразии диффеоморфизмами.}

Известно~\cite{MZ55}, что последнее утверждение влечет доказываемое следствие.\footnote{Вместо выделенного утверждения о группе $A_p$ можно использовать положительное решение
пятой проблемы Гильберта, см. ниже.}
%QED
\end{proof}

%\smallskip
Следствие \ref{hsc} доказано в 1946~г. Бохнером и Монтгомери~\cite[The\-o\-rem 2, p. 208]{MZ55}
более сложным образом.
Оно является \textit{гладкой} версией недоказанной гипотезы Гильберта"--~Смита
(формулировка которой получается из формулировки следствия заменой
%слова <<диффеоморфизмами>> на <<гомеоморфизмами>>
слов <<гладком>> и <<диффеоморфизмами>> на <<топологическом>> и <<гомеоморфизмами>>).
Доказательство \textit{липшицевой} версии см. в~\cite{RS97}.

Гипотеза Гильберта"--~Смита появилась при решении пятой проблемы Гильберта:
\textit{любая ли
локально евклидова топологическая группа является группой Ли?}
Сама пятая проблема Гильберта была положительно решена в 1952~г.~\cite{MZ55} (независимо Глизоном, а
также Монтгомери и Циппиным).

См. также \cite{AO07, HR08, OY03}.

%Мы не приводим доказательства указанных редукций.

 %\newpage
\subsection*{Доказательство\footnote{Это доказательство (оно появилось в~\cite{Sk07}) проще оригинального~\cite{RSS96, RS00} (хотя использует те же идеи).
Оно обобщает уже разобранное доказательство простейшего случая "--- утверждения \ref{p:diff-f}.}
теоремы \ref{thm:t5}}

%Мы советуем читателю разобрать это доказательство сначала для $M=\R^2$, чего
%достаточно для элементарных следствий (тогда конец этого абзаца можно
%пропустить).
То, что $N$ является дифференцируемым подмно\-го\-об\-ра\-зием в $M$ "---
локальное условие.
Поэтому можно считать, что $M=\R^m$.

Обозначим $|x|:=\sqrt{x_1^2+\dots+x_m^2}$ и
$$
B^{m,k}_l:=\{(x_1,\dots,x_m)\in\R^m\colon {-}l^2x_k<|x|<1/l\text{ и }
l^2|x_i|<|x|\text{ для }k<i\le m\}.$$
Тогда

$\bullet$
% ??? заменить запятую на точку с запятой?
% нет
$B^{m,m+1}_l$ есть проколотая внутренность $m$"=мерного шара радиуса $1/l$,

$\bullet$
$B^{m,k}_l$ есть открытый конус над $(1/l^2)$"=окрестностью $k$"=мерного
% ??? заменить ", и" на точку с запятой?
% нет
полушария в $(m-1)$"=мерной сфере радиуса $1/l$ для $1\le k\le m$, и

$\bullet$
$B^{m,0}_l=\emptyset$ для $l>m$.

Обозначим через $\mathrm{O}_m$ группу ортогональных преобразований пространства $\R^m$.

Возьмем наибольшее $k\ge0$, для которого
\begin{multline*}
(*)\quad \text{при любом $x\in N$ существуют такие $l>m$ и $A\in \mathrm{O}_m$, что}\\
(x+AB^{m,k}_l)\cap N=\emptyset.
\end{multline*}
(Неформально это значит, что $N$ является <<$(m-k)$"=мерно липшицевым>>.)
Такое $k$ существует, поскольку $(*)$~справедливо при $k=0$.
См. рис. \ref{ris.os}, где $m=2$ и $k=1$.

Если $k=m+1$, то $N$ состоит из изолированных точек и теорема доказана.
Поэтому будем считать, что $k\le m$.
% ??? вставить "мы" перед фиксируем?
% нет
Далее фиксируем $m$ и $k$ и опускаем их из обозначений конуса $B^{m,k}_l$.
Возьмем произвольную последовательность $\{A_l\}$, всюду плотную в $\mathrm{O}_m$.
 Обозначим
$$N_l:=\{x\in N\colon (x+A_lB_l)\cap N=\emptyset\}.$$
Ввиду условия~$(*)$ имеем $N=\bigcup\limits_{l=1}^\infty N_l$.
Нетрудно проверить, что \textit{$N_l$ замкнуто в $N$} (докажите или см. детали в
\cite[лемма 3.1]{RSS96}).
Значит, по теореме Бэра о категории некоторое $N_l$ содержит непустое
открытое в $N$ множество \cite[лемма 3.2]{RSS96}.

%Рассмотрим точку $a\in\R^m-N$.
%Так как $N$ локально компактно, то существует такое $r$, что шар с центром
%в $a$ и радиусом $r$ пересекает $N$ только в некоторых точках своей граничной
%сферы.
%Любая такая граничная точка лежит в $N_l$ для некоторого $l$.
%Поэтому из гладкой объемлемой однородности множества $N$ вытекает, что

Поэтому существуют точка $x\in N$ и замкнутый $m$"=мерный куб $I^m$ диаметра
меньше $1/l$ с центром в $x$, для которых $N':=N\cap I^m\subset N_l$.
%(рис. 2 и 3).
Тогда
$$
[(y+A_lB_l)\cup(y-A_lB_l)]\cap N'=\emptyset
\quad\mbox{при любом } y\in N'.
\leqno (**)
$$
Действительно, если $z\in(y-A_lB_l)\cap N'$, то
$y\in(z+A_lB_l)\cap N'\subset N_l$, что
невозможно.

Так как $N$ замкнуто, то можно считать, что $N'$ компактно.
Можно также считать, что некоторая $(m-k)$"=мерная грань $L$ куба $I^m$
перпендикулярна $k$"=мерной плоскости $A_l(\R^k\times\skew{-3}\vec0)$ ($L=I^m$ при $k=0$).
Обозначим через $p\colon I^m\to L$ ортогональную проекцию.

%Так как $N$ локально компактно, то $N'$ и вместе с ним $p(N')$ локально
%компактны.
%Поэтому существует шар в $L$, пересечение которого с $p(N')$ компактно.
%Уменьшим, если нужно, куб $I^m$, и будем считать, что этот шар совпадает с $L$.

%\smallskip
\textit{Первый случай: $p(N')$ содержит открытое в $L$ множество $U$.}
(Это заведомо так для $k=m$, когда все уже очевидно, и это заведомо не так для $k=0$.)
Из~$(**)$ следует, что $p$ является взаимно однозначным на $N'$, и что
обратное отображение $q\colon U\to N'$ липшицево.
Поэтому $q$ имеет точку дифференцируемости \cite[теорема 3.1.6]{Fe69}.
Тогда из дифференцируемой объемлемой однородности вытекает, что $q$
дифференцируемо в любой точке.
Поэтому условие из определения дифференцируемого подмногообразия выполнено в
одной точке множества $N$.
Тогда из дифференцируемой объемлемой однородности вытекает, что $N$ является
дифференцируемым подмно\-го\-об\-ра\-зием.

%\smallskip
\textit{Второй случай: $p(N')$ не содержит никакого открытого в $L$ множества.}
(Значит, $k<m$.)
Так как $p(N')$ не содержит открытого в $L$ множества, то существует точка
$a\in L-p(N')$, достаточно близкая к центру грани $L$ (точнее, расстояние от
которой до центра грани $L$ меньше четверти диаметра этой грани).
Так как $p(N')$ компактно, то расстояние от $a$ до $p(N')$ не равно нулю и
существует точка $z\in N'$, для которой $|a-p(z)|$ равно этому расстоянию.
Поскольку $a$ достаточно близко к центру грани $L$, то $p(z)$ лежит
\textit{внутри} грани $L$.
Тогда открытый $(m-k)$"=мерный шар $D\subset L$ с центром в $a$ и радиусом
$|a-p(z)|$ не пересекает $p(N')$.
Поэтому $p^{-1}(D)\cap N'=\emptyset$.
Ясно, что
$$(z+A_lB_l)\cup(z-A_lB_l)\cup p^{-1}(D)\supset z+A_lB^{m,k+1}_s\quad
\mbox{для некоторого }s.$$
Отсюда и из~$(**)$ следует, что $(z+A_sB^{m,k+1}_s)\cap N=\emptyset$.
Так как $N$ дифференцируемо объемлемо однородно, то при любом $x\in N$
существуют окрестности $Uz$ и $Ux$ точек $z$ и $x$ в $\R^m$ и диффеоморфизм
$h\colon Uz\to Ux$, переводящий $z$ в $x$ и $Uz\cap N$ в $Ux\cap N$.
Тогда по определению диффеоморфизма
$$h(Uz\cap(z+A_sB^{m,k+1}_s))\supset x+AB^{m,k+1}_u\quad\mbox{для некоторых }
A\in \mathrm{O}_m\mbox{ и }u>m.$$
Значит, $(*)$~выполнено с заменой $k$ на $k+1$.
Это противоречит максимальности числа $k$.\hfill$\Box$

\end{document}